\RequirePackage[l2tabu, orthodox]{nag}

\documentclass[12pt]{amsart}
\usepackage{fullpage,url,amssymb,enumerate,colonequals}
\usepackage[all]{xy} 
\usepackage{mathrsfs} 

\usepackage[OT2,T1]{fontenc}
\DeclareSymbolFont{cyrletters}{OT2}{wncyr}{m}{n}
\DeclareMathSymbol{\Sha}{\mathalpha}{cyrletters}{"58}

\usepackage{color}

\newcommand{\defi}[1]{\textsf{#1}} 

\newcommand{\C}{\mathbb{C}}
\newcommand{\F}{\mathbb{F}}
\newcommand{\G}{\mathbb{G}}

\newcommand{\Q}{\mathbb{Q}}
\newcommand{\R}{\mathbb{R}}
\newcommand{\Z}{\mathbb{Z}}

\newcommand{\Adeles}{\mathbf{A}}

\newcommand{\calH}{\mathcal{H}}

\newcommand{\EE}{\mathscr{E}}

\DeclareMathOperator{\alt}{alt}

\DeclareMathOperator{\Average}{Average}

\DeclareMathOperator{\coker}{coker}

\DeclareMathOperator{\height}{ht}
\DeclareMathOperator{\im}{im}

\DeclareMathOperator{\inv}{inv}

\DeclareMathOperator{\ord}{ord}

\DeclareMathOperator{\Prob}{Prob}

\DeclareMathOperator{\re}{Re}

\DeclareMathOperator{\rk}{rk}

\DeclareMathOperator{\Sel}{Sel}

\newcommand{\tors}{{\operatorname{tors}}}

\newcommand{\GL}{\operatorname{GL}}
\newcommand{\HH}{{\operatorname{H}}}

\newcommand{\M}{\operatorname{M}}

\newcommand{\injects}{\hookrightarrow}
\newcommand{\intersect}{\cap} 
\newcommand{\isom}{\simeq}

\newcommand{\tensor}{\otimes} 
\newcommand{\To}{\longrightarrow}
\newcommand{\Union}{\bigcup} 

\newtheorem{theorem}{Theorem}[section]

\newtheorem{proposition}[theorem]{Proposition}
\newtheorem{conjecture}[theorem]{Conjecture}

\theoremstyle{definition}

\newtheorem{question}[theorem]{Question}

\theoremstyle{remark}
\newtheorem{remark}[theorem]{Remark}

\makeatletter
\g@addto@macro\bfseries{\boldmath} 
\makeatother

\usepackage{microtype}

\usepackage[
	backref,
	pdfauthor={Bjorn Poonen}, 
]{hyperref}
\usepackage[alphabetic,backrefs,lite,nobysame]{amsrefs} 

\begin{document}

\title{Heuristics for the arithmetic of elliptic curves}
\author{Bjorn Poonen}
\thanks{This survey article has been submitted to the Proceedings of the 2018 ICM.  The writing of this article was supported in part by National Science Foundation grant DMS-1601946 and Simons Foundation grants \#402472 (to Bjorn Poonen) and \#550033.}
\address{Department of Mathematics, Massachusetts Institute of Technology, Cambridge, MA 02139-4307, USA}
\email{poonen@math.mit.edu}
\urladdr{\url{http://math.mit.edu/~poonen/}}
\date{November 30, 2017}

\begin{abstract}
This is an introduction to 
a probabilistic model for the arithmetic of elliptic curves, 
a model developed in a series of articles of the author
with Bhargava, Kane, Lenstra, Park, Rains, Voight, and Wood.
We discuss the theoretical evidence for the model,
and we make predictions about elliptic curves
based on corresponding theorems proved about the model.
In particular, the model suggests that all but finitely many
elliptic curves over $\Q$ have rank $\le 21$, which would imply
that the rank is uniformly bounded.
\end{abstract}

\maketitle

\section{Introduction}\label{S:introduction}

Let $E$ be an elliptic curve over $\Q$ 
(see \cite{SilvermanAEC2009} for basic definitions).
Let $E(\Q)$ be the set of rational points on $E$.
The group law on $E$ gives $E(\Q)$ the structure of an abelian group,
and Mordell proved that $E(\Q)$ is finitely generated \cite{Mordell1922};
let $\rk E(\Q)$ denote its rank.
The present survey article, 
based primarily on articles of the author
with Eric Rains~\cite{Poonen-Rains2012-selmer},
with Manjul Bhargava, Daniel M. Kane, Hendrik Lenstra, and Eric Rains~\cite{Bhargava-Kane-Lenstra-Poonen-Rains2015}, and 
with Jennifer Park, John Voight, and Melanie Matchett Wood~\cite{Park-Poonen-Voight-Wood-preprint}
is concerned with the following question:

\begin{question}
\label{Q:bounded rank}
Is $\rk E(\Q)$ bounded as $E$ varies over all elliptic curves over $\Q$?
\end{question}

Question~\ref{Q:bounded rank} was implicitly asked 
by Poincar\'e in 1901~\cite{Poincare1901}*{p.~173},
even before $E(\Q)$ was known to be finitely generated! 
Since then, many authors have put forth guesses,
and the folklore expectation has flip-flopped at least once;
see \cite{Poincare1950-Oeuvres5}*{p.~495, end of footnote~(${}^3$)}, 
\cite{Honda1960}*{p.~98}, 
\cite{Cassels1966-diophantine}*{p.~257}, 
Tate~\cite{Tate1974}*{p.~194}, 
\cite{Mestre1982}, 
\cite{Mestre1986}*{II.1.1 and II.1.2}, 
\cite{Brumer1992}*{Section~1}, 
\cite{Ulmer2002}*{Conjecture~10.5}, 
and \cite{Farmer-Gonek-Hughes2007}*{(5.20)}, 
or see \cite{Park-Poonen-Voight-Wood-preprint}*{Section~3.1} 
for a summary.

The present survey 
describes a probabilistic model for the arithmetic of elliptic curves, 
and presents theorems about the model that suggest that $\rk E(\Q) \le 21$ 
for all but finitely many elliptic curves $E$,
and hence that $\rk E(\Q)$ is bounded.
Ours is not the first heuristic for boundedness: 
there is one by Rubin and Silverberg for a family of quadratic twists 
\cite{Rubin-Silverberg2000}*{Remarks 5.1 and~5.2}, 
and another by Granville, 
discussed in \cite{Watkins-et-al2014}*{Section~11} 
and developed further in \cite{Watkins-discursus}.
Interestingly, the latter also suggests a bound of $21$.

Modeling ranks directly is challenging 
because there are few theorems about the distribution of ranks.
Also, although there exists extensive computational data that suggests answers 
to some questions 
(e.g., \cite{Balakrishnan-Ho-Kaplan-Spicer-Stein-Weigandt2016}),
it seems that far more data would be needed to suggest answers to others.
Therefore, instead of modeling ranks in isolation,
we model ranks, Selmer groups, and Shafarevich--Tate groups
simultaneously, so that we can calibrate and corroborate the model 
using a diverse collection of known results.

\section{The arithmetic of elliptic curves}

\subsection{Counting elliptic curves by height}

Every elliptic curve $E$ over $\Q$ is isomorphic to the projective closure 
of a unique curve $y^2=x^3+Ax+B$ in which $A$ and $B$
are integers with $4A^3+27B^2 \ne 0$ (the smoothness condition)
such that there is no prime $p$ such that $p^4|A$ and $p^6|B$.
Let $\EE$ be the set of elliptic curves of this form,
so $\EE$ contains one curve in each isomorphism class.
Define the \defi{height} of $E \in \EE$ by
\[
	\height E \colonequals \max(|4A^3|,|27B^2|).
\]
(This definition is specific to the ground field $\Q$,
but it has analogues over other number fields.)
Define 
\[
	\EE_{\le H} \colonequals \{E \in \EE : \height E \le H\}.
\]

Ignoring constant factors, 
we have about $H^{1/3}$ integers $A$ with $|4A^3| \le H$,
and $H^{1/2}$ integers $B$ with $|27B^2| \le H$.
A positive fraction of such pairs $(A,B)$
satisfy the smoothness condition and divisibility conditions,
so one should expect $\#\EE_{\le H}$ to be about $H^{1/3} H^{1/2} = H^{5/6}$.
In fact, an elementary sieve argument \cite{Brumer1992}*{Lemma~4.3} 
proves the following:

\begin{proposition}
\label{P:elliptic curves of bounded height}
We have 
\[
	\#\EE_{\le H} = (2^{4/3} 3^{-3/2} \zeta(10)^{-1} + o(1)) \; H^{5/6}
\]
as $H \to \infty$.
\end{proposition}

Define the \defi{density} of a subset $S \subseteq \EE$ as
\[
	\lim_{H\rightarrow\infty}
	\frac{\#(S\cap \EE_{\le H})}{\# \EE_{\le H}},
\]
if the limit exists.
For example, it is a theorem that 100\% of elliptic curves $E$ over $\Q$
have no nontrivial rational torsion points;
this statement is to be interpreted as saying that the density of the set
$S \colonequals \{E \in \EE: E(\Q)_{\tors} = 0\}$ is $1$
(even though there do exist $E$ with $E(\Q)_{\tors} \ne 0$).

\subsection{Elliptic curves over local fields}

Our model will be inspired by theorems and conjectures about
the arithmetic of elliptic curves over $\Q$.
But before studying elliptic curves over $\Q$,
we should thoroughly understand elliptic curves over local fields.

Let $\Q_v$ be the completion of $\Q$ at a place $v$.
There is a natural injective homomorphism 
$\inv \colon \HH^2(\Q_v,\G_m) \to \Q/\Z$
that is an isomorphism if $v$ is nonarchimedean.

Let $E$ be an elliptic curve over $\Q_v$.
Fix $n \ge 1$.
Taking Galois cohomology in the exact sequence
\[
	0 \To E[n] \To E \stackrel{n}\To E \To 0
\]
yields a homomorphism $E(\Q_v)/n E(\Q_v) \to \HH^1(\Q_v,E[n])$.
Let $W_v$ be its image.
If $v$ is a nonarchimedean place not dividing $n$
and $E$ has good reduction, 
then $W_v$ equals the subgroup of unramified classes in $\HH^1(\Q_v,E[n])$ 
\cite{Poonen-Rains2012-selmer}*{Proposition~4.13}.

The theory of the Heisenberg group \cite{MumfordTheta3}*{pp.~44--46}
yields an exact sequence
\[
	1 \To \G_m \To \calH \To E[n] \To 1,
\]
which induces a map of sets
\[
	q_v \colon \HH^1(\Q_v,E[n]) \To \HH^2(\Q_v,\G_m) \stackrel{\inv}\injects \Q/\Z.
\]
It turns out that $q_v$ is a quadratic form in the sense that
$q_v(x+y)-q_v(x)-q_v(y)$ is bi-additive~\cite{Zarhin1974}*{\S2}.
Moreover, $q_v|_{W_v}=0$ \cite{O'Neil2002}*{Proposition~2.3}.
In fact, using Tate local duality one can show that 
$W_v$ is a maximal isotropic subgroup of $\HH^1(\Q_v,E[n])$
with respect to $q_v$ \cite{Poonen-Rains2012-selmer}*{Proposition~4.11}.

\subsection{Selmer groups and Shafarevich--Tate groups}

Now let $E$ be an elliptic curve over $\Q$.
Let $\Adeles = \prod'_v (\Q_v,\Z_v)$ be the ad\`ele ring of $\Q$;
here $v$ ranges over nontrivial places of $\Q$,
Write $E(\Adeles)$ for $\prod_v E(\Q_v)/nE(\Q_v)$,
and write $\HH^1(\Adeles,E[n])$ for 
the restricted product $\prod'_v (\HH^1(\Q_v,E[n]),W_v)$.
We have a commutative diagram
\[
\xymatrix{
E(\Q)/nE(\Q) \ar[r] \ar[d] & \HH^1(\Q,E[n]) \ar[d]^{\beta} \\
E(\Adeles)/nE(\Adeles) \ar[r]^{\alpha} & \HH^1(\Adeles,E[n]). \\
}
\]
The \defi{$n$-Selmer group} is defined by
$\Sel_n E \colonequals \beta^{-1}(\im \alpha) \subseteq \HH^1(\Q,E[n])$.
(This is equivalent to the classical definition; we have only replaced
$\prod_v \HH^1(\Q_v,E[n])$ with a subgroup $\HH^1(\Adeles,E[n])$
into which $\alpha$ and $\beta$ map.)
The reason for defining $\Sel_n E$ is that 
it is a computable finite upper bound for (the image of) $E(\Q)/nE(\Q)$.
On the other hand, the \defi{Shafarevich--Tate group} is defined by
\[
	\Sha = \Sha(E) \colonequals \ker\left( \HH^1(\Q,E) \to \prod_v \HH^1(\Q_v,E) \right).
\]
It is a torsion abelian group with an alternating pairing 
\[
	[\;,\;] \colon \Sha \times \Sha \to \Q/\Z
\]
defined by Cassels.
Conjecturally, $\Sha$ is finite; in this case, $[\;,\;]$ is nondegenerate 
and $\#\Sha$ is a square \cite{Cassels1962-IV}.
The definitions easily yield an exact sequence
\begin{equation}
\label{E:Selmer-Sha}
	0 \To \frac{E(\Q)}{nE(\Q)} \To \Sel_n E \To \Sha[n] \To 0,
\end{equation}
so $\Sha[n]$ is measuring the difference between $\Sel_n E$
and the group $E(\Q)/nE(\Q)$ it is trying to approximate.

Each group in \eqref{E:Selmer-Sha} decomposes according to 
the factorization of $n$ into powers of distinct primes,
so let us restrict to the case in which $n=p^e$ for some prime $p$
and nonnegative integer $e$.
Taking the direct limit over $e$ yields an exact sequence 
\[
	0 \To E(\Q) \tensor \frac{\Q_p}{\Z_p} \To \Sel_{p^\infty} E 
	\To \Sha[p^\infty] \To 0
\]
of $\Z_p$-modules in which 
$\Sel_{p^\infty} E \colonequals \varinjlim \Sel_{p^e} E$
and $\Sha[p^\infty] \colonequals \Union_{e \ge 0} \Sha[p^e]$.
Moreover, one can show that if $E(\Q)[p]=0$ (as holds for 100\% of curves), 
then $\Sel_{p^e} E \to (\Sel_{p^\infty} E)[p^e]$ is an isomorphism 
(cf.~\cite{Bhargava-Kane-Lenstra-Poonen-Rains2015}*{Proposition~5.9(b)}),
so no information about the individual $p^e$-Selmer groups
has been lost in passing to the limit.

\subsection{The Selmer group as an intersection of maximal isotropic direct summands}
\label{S:Selmer group is intersection}

If $\xi = (\xi_v) \in \HH^1(\Adeles,E[n])$, 
then for all but finitely many $v$
we have $\xi_v \in W_v$ and hence $q_v(\xi_v)=0$,
so we may define $Q(\xi) \colonequals \sum_v q_v(\xi_v)$.
This defines a quadratic form $Q \colon \HH^1(\Adeles,E[n]) \to \Q/\Z$.

\begin{theorem}
\label{T:intersection of maximal isotropic subgroups}
\hfill
\begin{enumerate}[\upshape (a)]
\item Each of $\im \alpha$ and $\im \beta$ 
is a maximal isotropic subgroup of $\HH^1(\Adeles,E[n])$
with respect to $Q$ \cite{Poonen-Rains2012-selmer}*{Theorem~4.14(a)}.
\item\label{I:beta injective} 
If $n$ is prime 
or $G_\Q \to \GL_2(\Z/n\Z)$ is surjective
then $\beta$ is injective.
(See \cite{Poonen-Rains2012-selmer}*{Proposition~3.3(e)} 
and \cite{Bhargava-Kane-Lenstra-Poonen-Rains2015}*{Proposition~6.1}.)
\end{enumerate}
\end{theorem}

By definition, $\beta(\Sel_n E) = (\im \alpha) \intersect (\im \beta)$.
Thus, under either hypothesis in~\eqref{I:beta injective},
$\Sel_n E$ is isomorphic to 
an intersection of maximal isotropic subgroups of $\HH^1(\Adeles,E[n])$.

Moreover, $\im \alpha$ is a direct summand of $\HH^1(\Adeles,E[n])$ 
\cite{Bhargava-Kane-Lenstra-Poonen-Rains2015}*{Corollary~6.8}.
It is conjectured that $\im \beta$ is a direct summand as well,
at least for asymptotically 100\% of elliptic curves over $\Q$
\cite{Bhargava-Kane-Lenstra-Poonen-Rains2015}*{Conjecture~6.9}, 
and it could hold for all of them.

\subsection{The Birch and Swinnerton-Dyer conjecture}

See \cite{Wiles2006} 
for an introduction to the Birch and Swinnerton-Dyer
conjecture more detailed than what we present here,
and see \cite{Stein-Wuthrich2013}*{Section~8} 
for some more recent advances towards it.

Let $E \in \EE$.
To $E$ one can associate its \defi{$L$-function} $L(E,s)$,
a holomorphic function initially defined when $\re s$ is sufficiently large,
but known to extend to a holomorphic function on $\C$
(this is proved using the modularity of $E$).
Just as the Dirichlet analytic class number formula 
expresses the residue at $s=1$ 
of the Dedekind zeta function of a number field $k$ in terms of 
the arithmetic of $k$,
the Birch and Swinnerton-Dyer conjecture expresses the leading term
in the Taylor expansion of $L(E,s)$ around $s=1$
in terms of the arithmetic of $E$.
We will state it only in the case that $\rk E(\Q)=0$
since that is all that we will need.
In addition to the quantities previously associated to $E$, we need 
\begin{itemize}
\item the \defi{real period} $\Omega$, 
defined as the integral over $E(\R)$ of a certain $1$-form; and 
\item the \defi{Tamagawa number} $c_p$ for each finite prime $p$, 
a $p$-adic volume analogous to the real period.
\end{itemize}
Also define 
\[
	\Sha_0(E) \colonequals 
	\begin{cases}
		\#\Sha(E), &\textup{if $\rk E(\Q) = 0$;} \\
		0, &\textup{if $\rk E(\Q) > 0$.}
	\end{cases}
\]

\begin{conjecture}[The rank~0 part of the Birch and Swinnerton-Dyer conjecture]
\label{C:BSD}
If $E \in \EE$, then
\begin{equation}
\label{E:BSD}
	L(E,1) = \frac{\Sha_0 \; \Omega \; \prod_p c_p}{\#E(\Q)_{\tors}^2}.
\end{equation}
\end{conjecture}

\begin{remark}
In the case where the rank $r \colonequals \rk E(\Q)$ is greater than $0$, 
Conjecture~\ref{C:BSD} states only that $L(E,1)=0$, 
whereas the full Birch and Swinnerton-Dyer conjecture 
predicts that $\ord_{s=1} L(E,s)=r$ and predicts the leading coefficient
in the Taylor expansion of $L(E,s)$ at $s=1$.
\end{remark}

Let $H = \height E$.
Following Lang~\cite{Lang1983}
(see also \cite{Goldfeld-Szpiro1995}, \cite{deWeger1998}, \cite{Hindry2007}, 
\cite{Watkins2008-ExpMath}, and \cite{Hindry-Pacheco2016}),
we estimate the typical size of $\Sha_0$ 
by estimating all the other quantities in \eqref{E:BSD} as $H \to \infty$;
see \cite{Park-Poonen-Voight-Wood-preprint}*{Section~6} for details.
The upshot is that if we average over $E$ 
and ignore factors that are $H^{o(1)}$,
then \eqref{E:BSD} simplifies to $1 \sim \Sha_0 \, \Omega$
and we obtain $\Sha_0 \sim \Omega^{-1} \sim H^{1/12}$.
More precisely:
\begin{itemize}
\item $\prod_p c_p = H^{o(1)}$ \cite{deWeger1998}*{Theorem~3}, \cite{Hindry2007}*{Lemma~3.5}, \cite{Watkins2008-ExpMath}*{pp.~114--115}, \cite{Park-Poonen-Voight-Wood-preprint}*{Lemma~6.2.1}; 
\item $\#E(\Q)_{\tors} \le 16$ \cite{Mazur1977};
\item $\Omega = H^{-1/12 + o(1)}$ \cite{Hindry2007}*{Lemma 3.7}, \cite{Park-Poonen-Voight-Wood-preprint}*{Corollary~6.1.3}; and
\item
the Riemann hypothesis for $L(E,s)$ implies that $L(E,1) \le H^{o(1)}$ 
\cite{Iwaniec-Sarnak2000}*{p.~713}.
In fact, it is reasonable to expect 
$\underset{E \in \EE_{\le H}}{\Average}\; L(E,1) \asymp 1$.
(The symbol $\asymp$ means that the left side is bounded above and 
below by positive constants times the right side.)
\end{itemize}
Thus we expect
\begin{equation}
\label{E:average Sha0}
	\underset{E \in \EE_{\le H}}{\Average}\; \Sha_0(E) = H^{1/12+o(1)}
\end{equation}
as $H \to \infty$.

\section{Modeling elliptic curves over \texorpdfstring{$\Q$}{Q}}

\subsection{Modeling the \texorpdfstring{$p$}{p}-Selmer group}

According to Theorem~\ref{T:intersection of maximal isotropic subgroups},
$\Sel_p E$ is isomorphic to an intersection of maximal isotropic subspaces
in an infinite-dimensional quadratic space over $\F_p$.
So one might ask whether one could make sense of 
choosing maximal isotropic subspaces 
in an infinite-dimensional quadratic space at random,
so that one could intersect two of them to obtain 
a space whose distribution is conjectured to be that of $\Sel_p E$.
This can be done by equipping an infinite-dimensional quadratic space
with a locally compact topology \cite{Poonen-Rains2012-selmer}*{Section~2},
but the resulting distribution can be obtained more simply
by working within a $2n$-dimensional quadratic space 
and taking the limit as $n \to \infty$.
Now every nondegenerate $2n$-dimensional quadratic space 
with a maximal isotropic subspace is isomorphic to
the quadratic space $V_n=(\F_p^{2n},Q)$, where $Q$ is the quadratic form
\[
Q(x_1,\ldots,x_n,y_1,\ldots,y_n) \colonequals x_1 y_1 + \cdots + x_n y_n.
\]
Therefore we conjecture that the distribution of $\dim_{\F_p} \Sel_p E$
as $E$ varies over $\EE$ equals the limit as $n \to \infty$
of the distribution of the dimension of the intersection 
of two maximal isotropic subspaces in $V_n$
chosen uniformly at random from the finitely many possibilities.
The limit exists and can be computed explicitly;
this yields the formula on the right in the following:

\begin{conjecture}[\cite{Poonen-Rains2012-selmer}*{Conjecture~1.1}]
For each $s \ge 0$, 
the density of $\{E \in \EE: \dim_{\F_p} \Sel_p E = s\}$ equals
\begin{equation}
\label{E:Sel_p distribution}
	\prod_{j \ge 0} (1+p^{-j})^{-1} \prod_{j=1}^s \frac{p}{p^j-1}.
\end{equation}
\end{conjecture}

\begin{remark}
Let $E_d$ be the elliptic curve $dy^2=x^3-x$ over $\Q$.
Heath-Brown proved that the density of integers $d$
such that $\dim_{\F_2} \Sel_2 E_d - 2 = s$ equals
\[
	\prod_{j \ge 0} (1+2^{-j})^{-1} \prod_{j=1}^s \frac{2}{2^j-1},
\]
matching~\eqref{E:Sel_p distribution} for $p=2$ 
\cites{Heath-Brown1993,Heath-Brown1994}.
(The $-2$ is there to remove the ``causal'' contribution to $\dim \Sel_2 E_d$
coming from $E_d(\Q)[2]$.)
As we have explained, this result is a natural consequence
of the theory of Section~\ref{S:Selmer group is intersection},
but in fact Heath-Brown's result came first
and the theory was reverse engineered from it \cite{Poonen-Rains2012-selmer}!
Heath-Brown's result was extended by 
Swinnerton-Dyer \cite{Swinnerton-Dyer2008} and Kane \cite{Kane2013}
to the family of quadratic twists of any $E \in \EE$
with $E[2] \subseteq E(\Q)$ and no cyclic $4$-isogeny.
\end{remark}

\subsection{Modeling the \texorpdfstring{$p^e$}{p to the e}-Selmer group}

If $p$ is replaced by $p^e$, then we should replace $\F_p^{2n}$
by $V_n \colonequals ((\Z/p^e\Z)^{2n},Q)$.
But now there are different types of maximal isotropic subgroups
up to isomorphism.
For example, if $e=2$, then 
$(\Z/p^2\Z)^n \times \{0\}^n$ and $(p\Z/p^2\Z)^{2n}$ are both
maximal isotropic subgroups;
of these, only the first is a direct summand of $V_n$.
In what follows, we will use only direct summands,
for reasons to be explained at the end of this section.

\begin{conjecture}
\label{C:Sel p^e}
If we intersect two random maximal isotropic direct summands of 
$V_n \colonequals ((\Z/p^e\Z)^{2n},Q)$
and take the limit as $n \to \infty$ of the resulting distribution,
we obtain the distribution of $\Sel_{p^e} E$
as $E$ varies over $\EE$.
\end{conjecture}

For $m \ge 1$, let $\sigma(m)$ denote the sum of the positive divisors of $m$.
One can prove that the limit as $n \to \infty$
of the average size of the random intersection equals $\sigma(p^e)$,
and there is an analogous result 
for positive integers $m$ not of the form $p^e$ 
\cite{Bhargava-Kane-Lenstra-Poonen-Rains2015}*{Proposition~5.20}.
This suggests the following:

\begin{conjecture}[\cite{Poonen-Rains2012-selmer}*{Conjecture~1(b)}, \cite{Bhargava-Kane-Lenstra-Poonen-Rains2015}*{Section~5.7}, \cite{Bhargava-Shankar-4selmer}*{Conjecture~4}]
\label{C:average Sel_m}
For each positive integer $m$,
\[
	\underset{E \in \EE}\Average\; \#\Sel_m E = \sigma(m).
\]
(The average is interpreted as the limit as $H \to \infty$
of the average over $\EE_{\le H}$.)
\end{conjecture}

One could similarly compute the higher moments of the conjectural distribution;
see \cite{Poonen-Rains2012-selmer}*{Proposition~2.22(a)}
and \cite{Bhargava-Kane-Lenstra-Poonen-Rains2015}*{Section~5.5}.

There are several reasons why insisting upon direct summands 
in Conjecture~\ref{C:Sel p^e} seems right:
\begin{itemize}
\item
Conjecturally, both of the maximal isotropic subgroups arising in the 
arithmetic of the elliptic curve \emph{are} direct summands:
see the last paragraph of Section~\ref{S:Selmer group is intersection}.
\item
Requiring direct summands is essentially the only way to make the 
model for $\Sel_{p^e} E$ consistent with the model for $\Sel_p E$,
given that $\Sel_p E \isom (\Sel_{p^e} E)[p]$ for 100\% of curves 
\cite{Bhargava-Kane-Lenstra-Poonen-Rains2015}*{Remark~6.12}.
\item 
It leads to Conjecture~\ref{C:average Sel_m}, which has been proved
for $m \le 5$ \cites{Bhargava-Shankar-2selmer,Bhargava-Shankar-3selmer,Bhargava-Shankar-4selmer,Bhargava-Shankar-5selmer}.
\end{itemize}

\subsection{Modeling the \texorpdfstring{$p^\infty$}{p to the infty}-Selmer group and the Shafarevich--Tate group}

Choosing a maximal isotropic direct summand of $((\Z/p^e\Z)^{2n},Q)$
compatibly for all $e$ is equivalent to choosing a maximal isotropic
direct summand of the quadratic $\Z_p$-module $V_n \colonequals (\Z_p^{2n},Q)$.
This observation will lead us to a process that models
$\Sel_{p^e} E$ for all $e$ simultaneously,
or equivalently, that models $\Sel_{p^\infty} E$ directly.
To simplify notation, for any $\Z_p$-module $X$, 
let $X'$ denote $X \tensor \frac{\Q_p}{\Z_p}$;
if $X$ is a $\Z_p$-submodule of $V_n$,
then $X'$ is a $\Z_p$-submodule of $V_n'$.

Now choose maximal isotropic direct summands $Z$ and $W$ of $V_n$
with respect to the measure arising from taking
the inverse limit over $e$ of the uniform measure on the set of
maximal isotropic direct summands of $(\Z/p^e\Z)^{2n}$ 
\cite{Bhargava-Kane-Lenstra-Poonen-Rains2015}*{Sections 2 and~4};
then we conjecture that the limiting distribution of $Z' \intersect W'$
as $n \to \infty$ equals the distribution of $\Sel_{p^\infty} E$
as $E$ varies over $\EE$.
Again, the point is that this limiting distribution is compatible
with the the previously conjectured distribution for $\Sel_{p^e} E$
for each nonnegative integer $e$,
and the conjecture for $\Sel_{p^e} E$ 
was based on \emph{theorems} about Selmer groups of elliptic curves
(see Section~\ref{S:Selmer group is intersection}).

Even better, using the same ingredients, 
we can model $\rk E(\Q)$ and $\Sha[p^\infty]$ at the same time:
\begin{conjecture}[\cite{Bhargava-Kane-Lenstra-Poonen-Rains2015}*{Conjecture~1.3}]
\label{C:RST}
If we choose maximal isotropic direct summands $Z$ and $W$ 
of $(\Z_p^{2n},Q)$ at random as above,
and we define
\[
	R \colonequals (Z \intersect W)', \qquad 
	S \colonequals Z' \intersect W', \qquad
	T \colonequals S/R,
\]
then the limit as $n \to \infty$ of the distribution of the exact sequence
\[
	0 \To R \To S \To T \To 0
\]
equals the distribution of the sequence
\[
	0 \To E(\Q) \tensor \frac{\Q_p}{\Z_p} \To \Sel_{p^\infty} E 
	\To \Sha[p^\infty] \To 0
\]
as $E$ varies over $\EE$.
\end{conjecture}

There are several pieces of indirect evidence for 
the rank and $\Sha$ predictions of Conjecture~\ref{C:RST}: 
\begin{itemize}
\item Each of $R$ and $E(\Q) \tensor \frac{\Q_p}{\Z_p}$ is isomorphic to $(\Q_p/\Z_p)^r$
for some nonnegative integer $r$, called the \defi{$\Z_p$-corank}
of the module.
\item The $\Z_p$-corank of $R$ is $0$ or $1$, with probability $1/2$ each 
\cite{Bhargava-Kane-Lenstra-Poonen-Rains2015}*{Proposition~5.6}.
Likewise, a variant of a conjecture of Goldfeld 
(see \cite{Goldfeld1979}*{Conjecture~B} 
and \cites{Katz-Sarnak1999a,Katz-Sarnak1999b})
predicts that $\rk E(\Q)$ 
(which equals the $\Z_p$-corank of $E(\Q) \tensor \frac{\Q_p}{\Z_p}$) 
is $0$, $1$, $\ge 2$
with densities $1/2$, $1/2$, $0$, respectively.
\item The group $T$ is finite and carries a 
nondegenerate alternating pairing with values in $\Q_p/\Z_p$,
just as $\Sha[p^\infty]$ conjecturally does 
(the $p$-part of the Cassels pairing).
In particular, $\#T$ is a square.
\item 
Smith has proved a result analogous to Conjecture~\ref{C:RST} 
for the family of quadratic twists of any $E \in \EE$ 
with $E[2] \subseteq E(\Q)$ and no cyclic $4$-isogeny \cite{Smith-preprint}.
\end{itemize}

Further evidence is that there are in fact \emph{three} distributions
that have been conjectured to be the distribution of $\Sha[p^\infty]$
as $E$ varies over rank~$r$ elliptic curves,
and these three distributions coincide \cite{Bhargava-Kane-Lenstra-Poonen-Rains2015}*{Theorems 1.6(c) and~1.10(b)}.
This is so even in the cases with $r \ge 2$, 
which conjecturally occur with density $0$.
For a fixed nonnegative integer $r$, the three distributions are as follows:
\begin{enumerate}[\upshape 1.]
\item A distribution defined by Delaunay 
\cites{Delaunay2001,Delaunay2007,Delaunay-Jouhet2014a}, 
who adapted the Cohen--Lenstra heuristics for class groups 
\cite{Cohen-Lenstra1984}.
\item The limit as $n \to \infty$ of the distribution of 
$T \colonequals (Z' \intersect W')/(Z \intersect W)'$
when $(Z,W)$ is sampled from 
the set of pairs of maximal isotropic direct summands of $(\Z_p^{2n},Q)$ 
satisfying $\rk_{\Z_p}(Z \intersect W)=r$.
(This set of pairs is the set of $\Z_p$-points of a scheme of finite
type, so it carries a natural measure 
\cite{Bhargava-Kane-Lenstra-Poonen-Rains2015}*{Sections 2 and~4}.)
\item The limit as $n \to \infty$ through integers of the same parity as $r$
of the distribution of $(\coker A)_{\tors}$ when $A$ is sampled from the
space of matrices in $\M_n(\Z_p)$ satisfying $A^T=-A$ and $\rk_{\Z_p}(\ker A)=r$;
here $\ker A$ and $\coker A$ are defined by viewing $A$ as 
a $\Z_p$-linear homomorphism $\Z_p^n \to \Z_p^n$.
\end{enumerate}
The last of these is inspired by 
the theorem of Friedman and Washington \cite{Friedman-Washington1989}
that for each odd prime $p$, 
the limit as $n \to \infty$ of the distribution $\coker A$
for $A \in \M_n(\Z_p)$ chosen at random with respect to Haar measure
equals the distribution conjectured by Cohen and Lenstra to be
the distribution of the $p$-primary part of the class group 
of a varying imaginary quadratic field.

\subsection{Modeling the rank of an elliptic curve}

In the previous section, we saw in the third construction 
that conditioning on $\rk_{\Z_p}(\ker A)=r$
yields the conjectural distribution of $\Sha[p^\infty]$ for rank~$r$ curves.
The simplest possible explanation for this would be that 
sampling $A$ at random from 
$\M_n(\Z_p)_{\alt} \colonequals \{ A \in \M_n(\Z_p) : A^T=-A\}$
\emph{without} conditioning on $\rk_{\Z_p}(\ker A)$
caused $\rk_{\Z_p}(\ker A)$ to be distributed like 
the rank of an elliptic curve.

What is the distribution of $\rk_{\Z_p}(\ker A)$?
If $n$ is even,
then the locus in $\M_n(\Z_p)_{\alt}$
defined by $\det A=0$ is the set of $\Z_p$-points of 
a hypersurface, which has Haar measure~$0$,
so $\rk_{\Z_p}(\ker A) = 0$ with probability~$1$.
If $n$ is odd, however, 
then $\rk_{\Z_p}(\ker A)$ cannot be $0$,
because $n-\rk_{\Z_p}(\ker A)$ is the rank of $A$,
which is even for an alternating matrix.
For $n$ odd, it turns out that 
$\rk_{\Z_p}(\ker A)=1$ with probability~$1$.
If we imagine that $n$ was chosen large and with random parity,
then the result is that $\rk_{\Z_p}(\ker A)$ is $0$ or $1$,
with probability $1/2$ each.
This result agrees with the variant of Goldfeld's conjecture
mentioned above.
This model cannot, however, distinguish the relative frequencies
of curves of various ranks $\ge 2$, because in the model
the event $\rk_{\Z_p}(\ker A) \ge 2$ occurs with probability~$0$.

Therefore we propose a refined model
in which instead of sampling $A$ from $\M_n(\Z_p)_{\alt}$,
we sample $A$ from the set $\M_n(\Z)_{\alt,\le X}$ 
of matrices in $\M_n(\Z)_{\alt}$
with entries bounded by a number $X$ 
\emph{depending on the height $H$ of the elliptic curve being modeled},
tending to $\infty$ as $H \to \infty$.
This way, for elliptic curves of a given height $H$,
the model predicts a potentially positive but diminishing probability
of each rank $\ge 2$ (the probability that an integer point in a box 
lies on a certain subvariety), and we can quantify the rate at which this
probability tends to $0$ as $H \to \infty$ 
in order to count curves of height up to $H$ having each given rank.
In fact, we let $n$ grow with $H$ as well.

Here, more precisely, is the refined model.
To model an elliptic curve $E$ of height $H$,
using functions $\eta(H)$ and $X(H)$ to be specified later, 
we do the following:
\begin{enumerate}[\upshape 1.]
\item Choose $n$ to be an integer of size about $\eta(H)$
of random parity
(e.g., we could choose $n$ uniformly at random from
$\{\lceil \eta(H) \rceil,\lceil \eta(H) \rceil+1\}$).
\item Choose $A_E \in \M_n(\Z)_{\alt, \le X(H)}$ uniformly at random,
independently for each $E$.
\item Define random variables $\Sha_E' \colonequals (\coker A)_{\tors}$
and $\rk_E' \colonequals \rk_\Z(\ker A)$.
\end{enumerate}
Think of $\Sha_E'$ as the ``pseudo-Shafarevich--Tate group'' of $E$
and $\rk_E'$ as the ``pseudo-rank'' of $E$;
their behavior is intended to model the actual $\Sha$ and rank.

To complete the description of the model, 
we must specify the functions $\eta(H)$ and $X(H)$.
We do this by asking ``How large is $\Sha_0$ on average?'',
both in the model and in reality.
Recall from \eqref{E:average Sha0} that we expect
\begin{equation}
\label{E:average Sha0 redux}
	\underset{E \in \EE_{\le H}}{\Average}\; \Sha_0(E) = H^{1/12+o(1)}.
\end{equation}
Define 
\[
	\Sha_{E,0}' \colonequals 
	\begin{cases}
		\#\Sha_E', &\textup{if $\rk'_E = 0$;} \\
		0, &\textup{if $\rk'_E > 0$.}
	\end{cases}
\]
Using that the determinant of an $n \times n$ matrix 
is given by a polynomial of degree~$n$ in the entries,
we can prove that 
\begin{equation}
\label{E:average pseudo-Sha}
	\underset{E \in \EE_{\le H}}{\Average}\; \Sha_{E,0}' = X(H)^{\eta(H)(1+o(1))},
\end{equation}
assuming that $\eta(H)$ does not grow too quickly with $H$.
Comparing \eqref{E:average Sha0 redux} and \eqref{E:average pseudo-Sha}
suggests choosing $\eta(H)$ and $X(H)$ so that $X(H)^{\eta(H)} = H^{1/12+o(1)}$.
We assume this from now on.
It turns out that we will not need to know any more about $\eta(H)$
and $X(H)$ than this.

\subsection{Consequences of the model}

To see what distribution of ranks is predicted by the refined model,
we must calculate the distribution of ranks of alternating matrices
whose entries are integers with bounded absolute value;
the relevant result, 
whose proof is adapted from \cite{Eskin-Katznelson1995},
is the following:

\begin{theorem}[cf.~\cite{Park-Poonen-Voight-Wood-preprint}*{Theorem~9.1.1}]
\label{thm:EskinKatznelsonAlternating}
If $1 \leq r \leq n$ and $n-r$ is even, 
and $A \in \M_n(\Z)_{\alt, \le X}$ is chosen uniformly at random,
then 
\[
	\Prob(\rk(\ker A ) \ge r) \asymp_n (X^n)^{-(r-1)/2}.
\]
(The subscript $n$ on the symbol $\asymp$ means that the implied
constants depend on $n$.)
\end{theorem}

Theorem~\ref{thm:EskinKatznelsonAlternating} 
implies that for fixed $r \ge 1$ and $E \in \EE$ of height $H$,
\begin{equation}
\label{E:Prob of rank >=r}
	\Prob(\rk_E' \ge r) 
	= (X(H)^{\eta(H)})^{-(r-1)/2 + o(1)} 
	= H^{-(r-1)/24 + o(1)}.
\end{equation}
Using this, and the fact $\#\EE_{\le H} \asymp H^{5/6} = H^{20/24}$ 
(Proposition~\ref{P:elliptic curves of bounded height}),
we can now sum \eqref{E:Prob of rank >=r} over $E \in \EE_{\le H}$
to prove the following theorem about our model:

\begin{theorem}[\cite{Park-Poonen-Voight-Wood-preprint}*{Theorem~7.3.3}]
\label{T:rank 21}
The following hold with probability~$1$:
\begin{align*}
	\#\{ E \in \EE_{\le H} : \rk'_E = 0 \} &= H^{20/24+o(1)} \\
	\#\{ E \in \EE_{\le H} : \rk'_E = 1 \} &= H^{20/24+o(1)} \\
	\#\{ E \in \EE_{\le H} : \rk'_E \ge 2 \} &= H^{19/24+o(1)} \\
	\#\{ E \in \EE_{\le H} : \rk'_E \ge 3 \} &= H^{18/24+o(1)} \\
	&\vdots \\
	\#\{ E \in \EE_{\le H} : \rk'_E \ge 20 \} &= H^{1/24+o(1)} \\
	\#\{ E \in \EE_{\le H} : \rk'_E \ge 21 \} &\le H^{o(1)}, \\
	\#\{ E \in \EE : \rk'_E > 21 \} &\textup{ is finite}.
\end{align*}
\end{theorem}

This suggests the conjecture that the same statements hold for 
the \emph{actual} ranks of elliptic curves over $\Q$.
In particular, we conjecture that $\rk E(\Q)$ is uniformly bounded,
bounded by the maximum of the ranks of the 
conjecturally finitely many elliptic curves of rank $>21$.

\begin{remark}
Elkies has found infinitely many elliptic curves over $\Q$ of rank $\ge 19$,
and one of rank $\ge 28$; 
these have remained the records since 2006 \cite{Elkies2006}.
\end{remark}

\section{Generalizations}

\subsection{Elliptic curves over global fields}

What about elliptic curves over other global fields $K$?
Let $\EE_K$ be a set of representatives for the isomorphism classes
of elliptic curves over $K$.
Let $B_K \colonequals \limsup_{E \in \EE_K} \rk E(K)$.
For example, the conjecture suggested by Theorem~\ref{T:rank 21}
predicts that $20 \le B_{\Q} \le 21$.

\begin{theorem}[\cite{Tate-Shafarevich1967}, \cite{Ulmer2002}]
\label{T:global function field}
If $K$ is a global function field, then $B_K = \infty$.
\end{theorem} 

Even for number fields, $B_K$ can be arbitrarily large 
(but maybe still always finite):

\begin{theorem}[\cite{Park-Poonen-Voight-Wood-preprint}*{Theorem~12.4.2}]
\label{T:number fields}
There exist number fields $K$ of arbitrarily high degree 
such that $B_K \ge [K:\Q]$.
\end{theorem} 

Examples of number fields $K$ for which $B_K$ is large include fields 
in anticyclotomic towers and certain multiquadratic fields;
see \cite{Park-Poonen-Voight-Wood-preprint}*{Section~12.4}.

A naive adaptation of our heuristic 
(see \cite{Park-Poonen-Voight-Wood-preprint}*{Sections 12.2 and~12.3})
would suggest $20 \le B_K \le 21$ for every global field $K$, 
but Theorems \ref{T:global function field} and~\ref{T:number fields}
contradict this conclusion.
Our rationalization of this is that the elliptic curves of high rank 
in Theorems \ref{T:global function field} and~\ref{T:number fields} 
are special in that they are definable over a proper subfield of $K$,
and these special curves exhibit arithmetic phenomena 
that our model does not take into account.
To exclude these curves, let $\EE_K^\circ$ be the set of $E \in \EE_K$
such that $E$ is not a base change of a curve from a proper subfield,
and let $B_K^\circ \colonequals \limsup_{E \in \EE_K^\circ} \rk E(K)$.
It is possible that $B_K^\circ < \infty$ for every global field $K$.

\begin{remark}
On the other hand, it is not true that $B_K^\circ \le 21$ for all number fields,
as we now explain.
Shioda proved that $y^2=x^3+t^{360}+1$ 
has rank $68$ over $\C(t)$ \cite{Shioda1992}.
In fact, it has rank $68$ also over $K(t)$ for a suitable number field $K$.
For this $K$, specialization yields 
infinitely many elliptic curves in $\EE_K^{\circ}$ of rank $\ge 68$.
Thus $B_K^\circ \ge 68$.
See \cite{Park-Poonen-Voight-Wood-preprint}*{Remark~12.3.1} for details.
\end{remark}

\subsection{Abelian varieties}

\begin{question}
\label{Q:abelian varieties}
For abelian varieties $A$ over number fields $K$,
is there a bound on $\rk A(K)$ depending only on $\dim A$ and $[K:\Q]$?
\end{question}

By restriction of scalars, we can reduce to the case $K=\Q$
at the expense of increasing the dimension.
By ``Zarhin's trick'' that 
$A^4 \times (A^\vee)^4$ is principally polarized \cite{Zarhin1974-trick}, 
we can reduce to the case that $A$ is principally polarized,
again at the expense of increasing the dimension.
For fixed $g \ge 0$,
one can write down a family of projective varieties including
all $g$-dimensional principally polarized abelian varieties over $\Q$, 
probably with each isomorphism class represented infinitely many times.
We can assume that each abelian variety $A$ is defined by 
a system of homogeneous polynomials with integer coefficients,
in which the number of variables, the number of polynomials,
and their degrees are bounded in terms of $g$.
Define the height of $A$ to be the maximum of the absolute
values of the coefficients.
Then the number of $g$-dimensional principally polarized abelian varieties 
over $\Q$ of height $\le H$ is bounded by a polynomial in $H$.
If there is a model involving a pseudo-rank $\rk_A'$ 
whose probability of exceeding $r$ 
gets divided by at least a fixed fractional power of $H$
each time $r$ is incremented by $1$, 
as we had for elliptic curves,
then the pseudo-ranks are bounded with probability~$1$.
This might suggest a positive answer to Question~\ref{Q:abelian varieties},
though the evidence is much flimsier than in the case of elliptic curves.

\section*{Acknowledgments} 

I thank Nicolas Billerey, Serge Cantat, Andrew Granville, Eric Rains, Michael Stoll, and John Voight for comments.

\end{document}